\documentstyle[11pt,leqno]{article}
\newtheorem{theorem}{{\sc Theorem}}
\newcommand{\bt}{\begin{theorem}}
\newcommand{\et}{\end{theorem}}
\newcommand{\newsection}[1]{\setcounter{equation}{0} \setcounter{theorem}{0}
\section{#1}}

\newcommand{\NI}{\noindent}
\newcommand{\bea}{\begin{eqnarray}}
\newcommand{\eea}{\end{eqnarray}}

\def \spec#1 {\mathop{#1}}

\def \b #1 {\bf #1}

\newcommand {\CC}{\centerline}

\newcommand{\clf}{{\cal F}}

\newcommand{\ity}{\infty}
\newcommand{\raro}{\rightarrow}

\newcommand{\vsp}{\vskip 1em}

\newcommand{\be}{\begin{equation}}
\newcommand{\ee}{\end{equation}}
\newcommand{\ben}{\begin{eqnarray*}}
\newcommand{\een}{\end{eqnarray*}}

\begin{document}
\CC {\bf Maximal inequalities for bifractional Brownian motion}
\vsp
\CC{ B.L.S. Prakasa Rao}
\CC{CR Rao Advanced Institute of Mathematics, Statistics}
\CC{and Computer Science, Hyderabad 500046, India}
\vskip 1cm
\NI{\bf Abstract :} We derive some  maximal inequalities for the bifractional Brownian motion using comparison theorems for Gaussian processes.
\vsp
\NI{\bf Keywords and phrases:} bifractional Brownian motion; Fractional Brownian motion; Sudakov-Fernique  inequality; Maximal inequalities.
\vsp
\NI{\bf MSC 2010:} Primary 60G22.
\vsp
\newsection{Introduction} Fractional Brownian motion (fBm)  $W^H= \{W^H(t), t \geq 0\}$ has been used for modeling stochastic phenomena with long-range dependence. It is a centered Gaussian process with the covariance function
$$R_H(s,t)= \frac{1}{2}(t^{2H}+s^{2H}-|t-s|^{2H})$$
where $0<H<1$ and the constant $H$ is called the Hurst index. The case $H=1/2$ corresponds to the Brownian motion. For the properties of a fBm, see  Mishura (2008)  and Prakasa Rao (2010). Consider a centered Gaussian process $W^{H,K}= \{W_t^{H,K}, t\geq 0\}$ called  the {\it bifractional Brownian motion} (bifBm) with the covariance function
$$C_H(s,t)= \frac{1}{2^K}[(t^{2H}+s^{2H})^K-|s-t|^{2HK}]$$
where $0<H<1$ and  $0< K\leq 1.$ If $K=1$, then the bifractional Brownian motion reduces to  the fractional Brownian motion and if $K=1$ and $H=\frac{1}{2},$ then it reduces to the Brownian motion. A bifBm can be considered as a generalization of the fBm but its increments are not stationary. Russo and Tudor (2006) studied the properties of a bifbm $W^{H,K}$. Houdre and Villa (2003) and Tudor and Xiao (2007) discussed the following properties of a bifBM $W^{H,K}$ (cf. Tudor (2023)).

(1) $E(W_t^{H,K})=0, Var(W_t^{H,K})= t^{2HK}, t\geq 0.$ \\
(2) The process $W^{H,K}$ is self-similar with index $HK \in (0,1)$, that is, for every real $a>0,$
$$\{W_{at}^{H,K},t\geq 0\}=\{ a^{HK} W_t^{H,K},t \geq 0\}$$
in the sense that the processes, on both sides of the equality sign, have the same finite-dimensional distributions.\\
(3) The process $W^{H,K}$ is not Markov and it is not a semimartingale if $HK \neq \frac{1}{2}.$\\
(4) The sample paths of the process $W^{H,K}$ are Holder continuous of order $\delta$, for any $\delta<HK,$ and they are nowhere differentiable.\\
(5) The bifBm $W^{H,K}$ satisfies the inequalities
$$2^{-K}|t-s|^{2HK}\leq E[W_t^{H,K}-W_s^{H,K}]^2\leq 2^{1-K}|t-s|^{2HK}, t\geq 0, s \geq 0.$$
\vsp
The bifBm $W^{H,K}$ can be extended for $K\in (1,2)$ with $H \in (0,1)$ and $HK \in (0,1)$ (cf. Bardina and Es-Sebaiy (2011) and Lifshits and Volkava (2015)).

\newsection{\bf Maximal inequalities}

Some maximal inequalities, for a fractional Brownian motion with polynomial drift, are presented in Prakasa Rao (2013). Prakasa Rao (2014) gives an overview of maximal inequalities for a fractional Brownian motion. Our interest is to derive some maximal inequalities for the bifBm  as applications of the Sudakov-Fernique inequality (cf. Adler (1980) , Sudakov (1971, 1976); Fernique (1975,1997)) and their refinements by Vitale (2000). 
\vsp
For any process $X,$ defined on the underlying  probability space $(\Omega, \clf,P ),$ let $X^*$ denote the supremum process defined by
$$X_t^*=\sup_{0\leq s \leq t}|X_s|$$
whenever it is defined. Since the process $W^{H,K}$ is self-similar, it follows that
$$\{W_{at}^{H,K}, 0 \leq t \leq T\} \stackrel{\Delta} = \{a^{HK} W_t^{H,K}, 0 \leq t \leq T\}$$
for any $a>0$ and hence
$${W_{at}^{H,K}}^* \stackrel {\Delta}= a^{HK} {W_t^{H,K}}^*.$$
Here $\Delta$ indicates that the processes, on both sides  of the equations, have the same finite-dimensional distributions. We have the following result as a consequence of the self-similarity of the process $W^{H,K} $.
\vsp
\NI{\bf Theorem 2.1.} {\it For any $T>0$ and $p>0,$
$$E[({W^{H,K}}^*(T))^p]= K(H,p)T^{pHK}$$
where $K(H,p)=E[({W^{H,K}}^*(1))^p].$}
\vsp
The following result is due to Sudakov (1971,1976) and  Fernique (1975,1997). In the following discussion, we assume that the sample paths of the processes $X$ and $Y$ over the interval $[0,T]$ are continuous  almost surely. This will ensure that the suprema of different sets of random variables are measurable, if necessary, by choosing separable modifications of the processes involved (cf. Doob (1953)).
\vsp
\NI{\bf Theorem 2.2.} (Sudakov-Fernique inequality ) {\it Let the processes $X= \{X(t), 0\leq t \leq T\}$ and $Y=\{Y(t), 0\leq t \leq T\}$ be centered Gaussian processes such that
$$E|X_s-X_t|^2 \leq E|Y_s-Y_t|^2, 0\leq s,t \leq T.$$
Then}
\be
E[\sup_{0\leq t \leq T}X_t]\leq E[\sup_{0\leq t \leq T}Y_t].
\ee
\vsp
Choose the process $X$ as the bifractional Brownian motion $W^{H,K} $, the process $Y^{(1)}$ as the process $2^{(1-K)/2}W^{HK}$ and the process $Y^{(2)}$ as the process $2^{(-K)/2}W^{HK}$ where $W^{HK}$ is the fractional Brownian motion Hurst index $HK$ in Theorem 2.2. As an application of property (5) of a bifBm given in Section 1, it follows that 

$$ E|X_s-X_t|^2\leq E|Y^{(1)}_s-Y^{(1)}_t|^2, 0\leq t, s \leq T$$

and 
$$E|X_s-X_t|^2\geq E|Y^{(2)}_s-Y^{(2)}_t|^2, 0\leq t,s \leq T$$.

Applying the Sudakov-Fernique inequality, we obtain the following result.
\vsp
\NI{\bf Theorem 2.3.} {\it Let $\{W_t^{H,K}, 0\leq t \leq T\}$ be the bifractional Brownian motion and $\{W^{HK}(t), 0\leq t \leq T\}$ be the fractional Brownian motion with Hurst index $HK.$ Then 

\be
E[\sup_{0\leq t \leq T}W_t^{H,K}] \leq E[\sup_{0\leq t \leq T}Y^{(1)}_t].
\ee
and}
\be
E[\sup_{0\leq t \leq T}W_t^{H,K}] \geq E[\sup_{0\leq t \leq T}Y^{(2)}_t].
\ee

\NI{\bf Remarks:} Note that the random variables 
$$\sup_{0\leq t \leq r}W_t^{H,K} $$
and
$$ r^{HK} \sup_{0\leq t \leq 1}W_t^{H,K}$$
have the same probability distribution by the self-similarity of the process $\{W_t^{H,K}, t \geq 0\}.$  Hence
$$P(\sup_{0\leq t \leq r} W_t^{H,K}\geq u)=P(\sup_{0\leq t \leq 1} W_t^{H,K}\geq u r^{-HK})$$ 
for any $r>0.$ From the symmetry property of the bifractional Brownian motion, it follows that
\bea
P(\sup_{0\leq t \leq T}W_t^{H,K}\geq u) &= & P(-\sup_{0\leq t \leq T}W_t^{H,K}\leq -u)\\\nonumber
&=& P(\inf_{0\leq t \leq T}(-W_t^{H,K})\leq -u)\\\nonumber
&=& P(\inf_{0\leq t \leq T}W_t^{H,K}\leq -u).\\\nonumber
\eea
Hence, for any $u>0,$
\bea
P(\sup_{0\leq t \leq T}|W_t^{H,K}(t)|\geq u) &\leq & P(\sup_{0\leq t \leq T}W_t^{H,K}\geq u)\\\nonumber
&&\;\;\;\;+P(\inf_{0\leq t \leq T}W_t^{H,K}\leq -u)\\\nonumber
&=& 2 P(\sup_{0\leq t \leq T}W_t^{H,K}\geq u).\\\nonumber
\eea

The following result is a variant of the Sudakov-Fernique inequality (cf. Alexander (1985)).

\NI{\bf Theorem 2.4.} {\it Suppose that $\{X_t, 0\leq t \leq T\}$ and $\{Y_t, 0\leq t \leq T\}$ are two centered Gaussian processes. Further suppose that
\be
E(X_t-X_s)^2 \leq E(Y_t-Y_s)^2, 0\leq t,s \leq T.
\ee
Then, for any non-decreasing convex function $g(.): R_+ \raro R,$ }
\be
E[g(\sup_{0\leq t,s \leq T}(X_t-X_s))]\leq E[g(\sup_{0\leq t,s \leq T}(Y_t-Y_s))]
\ee
and
\be
E[\sup_{0\leq t \leq T}X_t]\leq E[\sup_{0\leq t \leq T}Y_t].
\ee
\vsp
Applying Theorem 2.4 to the processes $X=W^{H,K}, Y^{(1)}=2^{(1-K)/2}W^{HK}$, and $Y^{(2)}=2^{(-K)/2}W^{HK}$, we obtain the following results.
\vsp
\NI{\bf Theorem 2.5.} {\it Let  $\{W_t^{H,K}, 0\leq t \leq T\}$ be the bifractional Brownian motion and $\{Y_t^{(1)}, 0\leq t \leq T\}$ be the fractional Brownian motion as defined earlier. Then, for any non-decreasing convex function $g(.): R_+ \raro R,$ 
\be
E[g(\sup_{0 \leq t,s \leq T}(W_t^{H,K}-W_s^{H,K})]\leq E[g(\sup_{0\leq t,s \leq T}(Y_t^{(1)}-Y_s^{(1)})]
\ee
and}
\be
E[\sup_{0\leq t \leq T}W_t^{H,K}]\leq E[\sup_{0\leq t \leq T}Y_t^{(1)}].
\ee
\vsp
\NI{\bf Theorem 2.6.} {\it Let  $\{W_t^{H,K}, 0\leq t \leq T\}$ be the bifractional Brownian motion and $\{Y_t^{(2)}(t), 0\leq t \leq T\}$ be the fractional Brownian motion as defined earlier. Then, for any non-decreasing convex function $g(.): R_+ \raro R,$ 
\be
E[g(\sup_{0\leq t,s \leq T}(W_t^{H,K}-W_s^{H,K})]\geq E[g(\sup_{0\leq t,s \leq T}(Y_t^{(2)}-Y_s^{(2)})]
\ee
and}
\be
E[\sup_{0\leq t \leq T}W_t^{H,K}]\geq E[\sup_{0\leq t \leq T}Y_t^{(2)}].
\ee
\vsp
Vitale (1996) extended Theorem 2.4  using the fact that any positive constant can be considered as nearly the supremum of a centered Gaussian process in the following sense.
\vsp
\NI{\bf Lemma 2.7.} {\it For any $c\geq 0,$ there is a sequence of mean zero Gaussian random variables such that (i) $c\leq \sup_{j\geq 1}W_j<\ity$ and (ii) $\sup_{j\geq n}W_j $ decreases to $c$ as $n \raro \ity.$} 
\vsp
For proof of Lemma 2.7, see Vitale (1996). The following theorem is due to Vitale (2000).
\vsp
\NI{\bf Theorem 2.8.}  {\it Suppose that $\{X_t, 0\leq t \leq T\}$ and $\{Y_t, 0\leq t \leq T\}$ are two centered Gaussian processes. Further suppose that
\be
E(X_t-X_s)^2 \leq E(Y_t-Y_s)^2, 0\leq t,s \leq T.
\ee
Let $m(t)$ be any continuous  function. Then 
\be
E[\sup_{0\leq t \leq T}(X_t+m(t))]\leq E[\sup_{0\leq t \leq T}(Y_t+m(t))].
\ee
Furthermore, for any constant $c$ and fixed $s\in [0,T],$
\be
E[\sup_{0\leq t \leq T}\{X_t+m(t)-X_s,c\}]\leq E[\sup_{0\leq t \leq T}\{Y_t+m(t)-Y_s,c\}].
\ee
In addition, suppose that $g:R\raro R$ is a non-decreasing convex function such that, for any fixed $s \in [0,T],$
$$E[g(\sup_{0\leq t \leq T}(Y_t+m(t)-Y_s))]<\ity.$$
Then }
\be
E[g(\sup_{0\leq t \leq T}(X_t+m(t)-X_s))]\leq E[g(\sup_{0\leq t \leq T}(Y_t+m(t)-Y_s))].
\ee
\vsp
As a consequence of  Theorem 2.8, we obtain the following result for a bifractional Brownian motion. 
\vsp
\NI{\bf Theorem 2.9.} {\it Let $W^{H,K}$ be a  bifractional Brownian motion with Hurst indices $H,K$ and define the processes $Y^{(i)}, i=1,2$ as given earlier.  Let $m(t)$ be any continuous function. Then 
\be
E[\sup_{0\leq t \leq T}(W_t^{H,K}+m(t))]\leq E[\sup_{0\leq t \leq T}(Y_t^{(1)}+m(t))] 
\ee
and
\be
E[\sup_{0\leq t \leq T}(W_t^{H,K}+m(t))]\geq E[\sup_{0\leq t \geq T}(Y_t^{(2)}+m(t))].
\ee
Furthermore, for any constant $c$ and fixed $s\in [0,T],$
\bea
\;\;\;\\\nonumber
E[\sup_{0\leq t \leq T}\{W_t^{H,K}+m(t)-W_s^{H,K},c\}] &\leq & E[\sup_{0\leq t \leq T}\{Y_t^{(1)}+m(t)-Y_s^{(1)},c\}]\\\nonumber
\eea
and
\bea
\;\;\;\\\nonumber
E[\sup_{0\leq t \leq T}\{W_t^{H,K}+m(t)-W_s^{H,K},c\}] &\geq & E[\sup_{0\leq t \leq T}\{Y_t^{(2)}+m(t)-Y_s^{(2)},c\}]\\\nonumber
\eea
In addition, suppose that $g(.):R\raro R$ is a non-decreasing convex function and fix $s \in [0,T].$. Then
\bea
\;\;\;\\\nonumber
E[g(\sup_{0\leq t \leq T}(W_t^{H,K}+m(t)-W_s^{H,K})] &\leq &E[g(\sup_{0\leq t\leq T}(Y_t^{(1)}+m(t)-Y_s^{(1)})] 
\\\nonumber
\eea
and
\bea
\;\;\;\;\\\nonumber
E[g(\sup_{0\leq t \leq T}(W_t^{H,K}+m(t)-W_s^{H,K})] &\geq & E[g(\sup_{0\leq t \leq T}(Y_t^{(2)}+m(t)-Y_s^{(2)})] 
\\\nonumber
\eea
provided that the expectations in the inequalities stated above exist.}
\vsp
Note that Theorem 2.9 is a consequence  of Theorem 2.8 since the condition (2.13) holds for the bifBm $W^{H,K}$ and the processes$Y^{(i)}, i=1,2$ by the property (5) in Section 1. Vitale (1996,2000) has proved a modified version of Theorem 2.4. We now state this result.
\vsp
\NI{\bf Theorem 2.10.}  {\it Suppose that $\{X_t, 0\leq t \leq T\}$ and $\{Y_t, 0\leq t \leq T\}$ are two centered Gaussian processes. Further suppose that
\be
E(X_t-X_s)^2 \leq E(Y_t-Y_s)^2, 0\leq t,s \leq T
\ee
and
\be
E(X_t^2) \leq E(Y_t^2), 0\leq t \leq T.
\ee
Let $m(t)$ be any continuous function and $c$ be an arbitrary constant. Then 
\be
E[\sup_{0\leq t \leq T}\{X_t+m(t),c\}]\leq E[\sup_{0\leq t \leq T}\{Y_t+m(t),c\}].
\ee
In addition, suppose that $g(.):R\raro R$ is a non-decreasing convex function.Then 
\be
E[g(\sup_{0\leq t\leq T}(X_t+m(t))]\leq E[g(\sup_{0\leq t \leq T}(Y_t+m(t))
\ee
and, for every $0\leq t \leq T,$}
\be
\int_t^\ity P(\sup_{0\leq s \leq T}X_s >u)du \leq \int_t^\ity P(\sup_{0\leq s \leq T}Y_s >u)du.
\ee
\vsp
Let $X=W^{H,K}$ and $Y^{(i)}, i=1,2$ be the processes as stated above . Note that $E(X_t^2)= t^{2HK}$.  An application of Theorem 2.10 leads to the following results for a bifractional Brownian motion.
\vsp
\NI{\bf Theorem 2.11.} {\it Let $X=W^{H,K}$ and $Y^{(1)}$ be the process as defined earlier.Then the following inequalities hold. Let $m(t)$ be any continuous function and $c$ be an arbitrary constant. Then 
\be
E[\sup_{0\leq t \leq T}(W_t^{H,K}+m(t),c)]\leq E[\sup_{0\leq t \leq T}(Y_t^{(1)}+m(t),c)]
\ee
In addition, suppose that $g:R\raro R$ is a non-decreasing convex function.Then 
\be
E[g(\sup_{0\leq t \leq T}(W_t^{H,K}+m(t))]\leq E[g(\sup_{0\leq t \leq T}(Y_t^{(1)}+m(t))].
\ee
and, for any $0\leq t \leq T,$
\be
\int_t^\ity P(\sup_{0\leq s \leq T}W_s^{H,K} >u)du \leq \int_t^\ity P(\sup_{0\leq s \leq T}Y_s^{(1)} >u)du.
\ee
or equivalently, for any $0\leq t \leq T,$}
\be
E[\sup_{0\leq s \leq T}W_s^{H,K}-t]_+\leq E[\sup_{0\leq s \leq T}Y_s^{(1)}-t]_+.
\ee
\vsp
\NI{\bf Theorem 2.12.} {\it Let $X=W^{H,K}$ and $Y^{(2)}$ be the process as defined earlier. Then the following inequalities hold.
Let $m(t)$ be any continuous  function and $c$ be an arbitrary constant. Then 
\be
E[\sup_{0\leq t \leq T}\{W_t^{H,K}+m(t),c\}]\geq E[\sup_{0\leq t \leq T}\{Y_t^{(2)}+m(t),c\}]
\ee
In addition, suppose that $g(.):R\raro R$ is a non-decreasing convex function.Then 
\be
E[g(\sup_{0\leq t \leq T}(W_t^{H,K}+m(t))]\geq E[g(\sup_{0\leq t\leq T}(Y_t^{(2)}+m(t))].
\ee
and, for any $0\leq t \leq T,$
\be
\int_t^\ity P(\sup_{0\leq s \leq T}W_s^{H,K} >u)du \geq \int_t^\ity P(\sup_{0\leq s \leq T}Y_s^{(2)} >u)du.
\ee
or equivalently for any $0\leq t \leq T,$}
\be
E[\sup_{0\leq s \leq T}W_s^{H,K}-t]_+\geq E[\sup_{0\leq s \leq T}Y_s^{(2)}-t]_+.
\ee
\vsp
\NI{\bf Remarks:} Vitale (1996) proved analogues of theorems stated above for denumerable sequences of Gaussian random variables. The results continue to hold for Gaussian processes with continuous sample paths almost surely which can be assumed to be separable through standard modifications (cf. Doob (1953)) and the results stated above follow from the monotone convergence theorem. Observe that the bifractional Brownian motion and the fractional Brownian motion have seperable modifications and the seperable versions, over a finite interval, have continuous sample paths almost surely.
\vsp
\vsp
Vitale (1996) obtained the following result giving an exponential inequality for centered Gaussian processes using the concept of Wills functional (cf. Wills (1973)) via geometric methods. 
\vsp
\NI{\bf Theorem 2.13.} {\it Suppose that $\{X_t,0\leq t \leq T\}$ is a centered Gaussian process and bounded in the sense that $E[\sup_{0\leq t \leq T}|X_t|]<\ity.$ Let $\sigma_t^2= E(X_t^2).$ Then
\be
E(\exp(\sup_{0\leq t\leq T}[X_t-\frac{1}{2}\sigma_t^2]))\leq \exp(E[\sup_{0\leq t \leq T}X_t])
\ee
and, for any $a>0,$
\be
P[\sup_{0\leq t \leq T}X_t-E(\sup_{0\leq t \leq T}X_t) \geq a] \leq \exp[-(1/2)(a^2/\sigma^2)]
\ee
where $\sigma^2= \sup_{0\leq t \leq T}\sigma_t^2.$}
\vsp
As special cases of Theorem 2.13, we have the following results for a bifractional Brownian motion $\{W^{H,K}_t,0\leq t \leq T\}$ observing that it is a centered Gaussian processes with $E[W_t^{H,K}]^2= t^{2HK}.$
\vsp
\NI{\bf Theorem 2.14.} {\it Suppose that $\{W_t^{H,K},0\leq t \leq T\}$ is a bifractional Brownian motion. Then
\be
E(\exp(\sup_{0\leq t\leq T}[W_t^{H,K}-\frac{1}{2}t^{2HK}]))\leq \exp(E[\sup_{0\leq t \leq T}W_t^{H,K}])
\ee
and, for any $a>0,$}
\be
P[\sup_{0\leq t \leq T}W^{H,K}_t-E(\sup_{0\leq t \leq T}W^{H,K}_t) \geq a] \leq \exp[-(1/2)(a^2/T^{2HK})].
\ee
\vsp
\NI{\bf Acknowledgment:} This work was supported by the Indian National Science Academy (INSA) under the scheme  ``INSA Honarary Scientist" at the CR RAO Advanced Institute of Mathematics, Statistics and Computer Science, Hyderabad 500046, India.

\NI{\bf References:}
\begin{description}

\item Adler, R.J. {\it  An Introduction to Continuity, Extrema, and Related Topics for General Gaussian Processes}, Institute of Mathematical Statistics, Lecture Notes-Monograph Series, Vol. 2, Hayward, California; 1990.

\item  Alexander, R. Lipschitzian mappings and total mean curvature of polyhedral surfaces I,  Trans. Amer. Math. Soc. 1985; 288: 661-678.

\item Bardina, X. and Es-Sebaiy, K. An extension of bifractional Brownian motion, Communications on Stochastic Analysis. 2011: 5: 333-340.

\item Doob, J.L. Stochastic Processes, Wiley, New York; 1953.

\item Fernique, X. Regularite des trajectoires des fonctions aleatoires Gauussiennes, Lecture Notes in Mathematics. 1975; 480: 1-96, Springer, Berlin.

\item Fernique, X. Fonctions aleatoires gaussiennes vecteurs aleatoires Gaussiens, CRM, Montreal; 1997.

\item Houdre, C. and Villa, J. An example of infinite dimensional quasi-helix, Contemp. Math. 2003; 336:195-201. 

\item Lipshits, M. and Volkova, K. Bifractional Brownian motion: Existence and border cases, arXiv:1502.02217, 2015. 

\item Mishura, Y. {\it Stochastic Calculus for Fractional Brownian Motion and Related Processes}, Lecture Notes in Mathematics No. 1929,  Springer-Verlag, Berlin; 2008.

\item Prakasa Rao, B.L.S.  {\it Statistical Inference for Fractional Diffusion Processes}, Wiley, London, 2010.

\item Prakasa Rao, B.L.S. Some maximal inequalities for fractional Brownian motion with polynomial drift, Stochastic Anal. and Appl. 2013; 31: 758-799.

\item Prakasa Rao, B.l.S. Maximal inequalities for fractional Brownian motion: an overview,  Stochastic Anal. and Appl. 2014; 32: 450-479.

\item Russo, F. and Tudor, C. On the bifractional Brownian motion, Stoch. Process. Their Appl. 2006; 116, 830-856.

\item Sudakov, V.N. Gaussian random processes and measures of solid angles in Hilbert spaces,  Soviet Math. Dokl. 1971; 12: 412-415.

\item Sudakov, V.N.  {\it Geometric Problems in the Theory of Infinite-Dimensional Probability Distributions}, Proc. Steklov Inst. Math. 1976; Vol.2, Amer. Math. Soc., Providence, Rhode Island.

\item Tudor,C. {\it Stochastic Partial Differential Equations with Additive Gaussian Noise}, World Scientific, Singapore, 2023.

\item Tudor, C. and Xiao, Y. Sample path properties of bifractional Brownian motion, Bernoulli. 2007; 13: 1023-1052.

\item Vitale, R.  Covariance identities for normal random variables via convex polytopes,  Statist. Probab. Lett. 1996; 30: 363-368.

\item Vitale, R. The Wills functional and Gaussian processes,  Ann. Probab. 1996; 24: 2172-2178.

\item Vitale, R. Some comparisons for Gaussian Processes, Proc. Amer. Math. Soc. 2000; 128: 3043-3046.

\item Wills, J.M. Zur Gitterpunktanzahl konvexer Mengen, Elem. Math. 1973; 28: 57-63.

\end{description}

\end{document}